\documentstyle[amssymb,amstex,graphicx,graphics,12pt,a4wide]{article}
\newtheorem{theorem}{Theorem}

\newtheorem{proposition}{Proposition}
\newtheorem{corollary}{Corollary}
\newtheorem{remark}{Remark}

\newcommand{\Faces}{\operatorname{Faces}}

\newcommand{\MaxFaces}{\operatorname{MaxPolytopes}}
\newcommand{\MaxFacets}{\operatorname{MaxPolytopes}}

\newcommand{\Cells}{\operatorname{Cells}}

\newcommand{\TypicalCell}{\operatorname{TypicalCell}}

\newcommand{\Vol}{\operatorname{Vol}}

\newcommand{\Cov}{\operatorname{Cov}}

\newcommand{\Var}{\operatorname{Var}}

\newcommand{\PHT}{\operatorname{PHT}}

\begin{document}

\title{Intrinsic Volumes of the Maximal Polytope Process in Higher Dimensional STIT Tessellations}
\author{Tomasz Schreiber (\textasteriskcentered 1975 -- \dag 2010)\\
Nicolaus Copernicus University, Toru\'n, Poland\\
Christoph Th\"ale\footnote{Current Address: University of Osnabr\"uck, Germany.\newline \texttt{e-mail: christoph.thaele[at]uni-osnabrueck.de}}\\
University of Fribourg, Fribourg, Switzerland}
\date{}
\maketitle

\begin{abstract}
Stationary and isotropic iteration stable random tessellations are considered, which can be constructed by a random process of cell division. The collection of maximal polytopes at a fixed time $t$ within a convex window $W\subset{\Bbb R}^d$ is regarded and formulas for mean values, variances, as well as a characterization of certain covariance measures are proved. The focus is on the case $d\geq 3$, which is different from the planar one, treated separately in \cite{ST2}. Moreover, a multivariate limit theorem for the vector of suitably rescaled intrinsic volumes is established, leading in each component -- in sharp contrast to the situation in the plane -- to a non-Gaussian limit.
\end{abstract}
\begin{flushleft}\footnotesize
\end{flushleft}

\section{Introduction}\label{secIntro}

Random tessellations {\bf st}able with respect to {\bf it}erations (STIT) have recently been introduced in \cite{MNW} and \cite{NW05} as a new model for random tessellation in ${\Bbb R}^d$ and have quickly attracted considerable interest in modern stochastic geometry as well fitting the growing demand for non-trivial and flexible but 
mathematically tractable tessellation models. The STIT tessellations may be interpreted as outcome of a random cell division process, which makes them very attractive for applications, see for example \cite{NMOW}. Other potential applications include mathematical modelling of systems of cracks, joints or fissures in rock or the so-called craquel\'ee of thin layers. A general approach to random cell division processes has recently appeared in \cite{Cowan} and the construction of STIT tessellations can roughly be described as follows. At first, we fix a compact and convex window $W\subset{\Bbb R}^d$ in which the construction is carried out. For simplicity we assume $W$ to be a polytope, as this implies that the resulting tessellation has polytopal cells with probability one. Next, an exponentially random life time is assigned to $W$, whereby the parameter of the distribution is given as a constant multiple of the integral-geometric mean width of $W$, see Section \ref{secSTIT} for details. Upon expiry of this life time a $W$ hitting uniform random hyperplane is chosen, is introduced in $W$ and divides the window into two polytopal sub-cells $W^+$ and $W^-$. The construction continues now recursively and independently in both of these sub-cells, where the newly introduced hyperplanes are always chopped-off by the boundaries of their mother-cells. The whole construction is continued until some deterministic time threshold $t$ is reached. Regarded in time, the construction can by interpreted as a pure-jump Markov process on the space of tessellations of the window $W$.\\ Our construction shares some common features with random fragmentation processes or branching Markov chains in the sense of \cite{Bertoin}. The cells of the resulting tessellation within $W$ can be regarded as particles in a suitable Polish space and the dynamics of the particles is non-interacting in a sense that different particles (cells of the tessellation) have independent evolutions, which is indeed the case in our construction. Moreover, whenever a particle dies it is replaced by exactly two new particles, namely the two newly generated sub-cells. In addition, the life times of the particles are exponentially distributed, as assumed in \cite{Bertoin}, but in general not independent in contrast to the fragmentation theory.\\ In the recent paper \cite{ST} the authors have introduced a new technique for studying the geometric properties of STIT tessellations based on martingales and the general theory of martingale problems for pure-jump type Markov processes. In particular, with these new developments, the variance of the total surface area of a stationary and isotropic iteration stable random tessellation $Y(t,W)$ in $W\subset{\Bbb R}^d$ (a random STIT tessellation) has been determined by integral-geometric means and the corresponding central limit theory has been established. Strikingly, as already signalled by results in the special
two dimensional case \cite{ST2}, and as confirmed by the results of the present paper, it turns out that
in a certain rather strong sense the asymptotic behaviour of the surface area process dominates and fully determines
the asymptotic geometry of the STIT tessellation. In particular the surface area variance is \textit{the} basic
second-order parameter of the tessellation and all second-order characteristics of functionals considered
in our work can be reduced to it, likewise non-trivial functional limits in law exist for the STIT surface area
process whereas the limits in law for other natural related processes arise as deterministic functionals
of the corresponding surface process. Another crucial phenomenon arising in this context, as first noted
in \cite{ST}, is that the asymptotic theories for dimension $d=2$ and $d>2$ differ strongly in many important
aspects. We have studied the planar case $d=2$ in the recent separate paper \cite{ST2} and thus only consider $d>2$ in the present one.\\ The purpose of this paper is to establish a second-order and limit theory for integral-geometric characteristics of stationary and isotropic STIT tessellations in dimensions higher than two. More precisely, the characteristics studied in this work are the cumulative intrinsic volumes of all orders for the collections of so-called maximal polytopes of $Y(t,W)$. These are in codimension $1$ the cell-separating facets introduced during the random cell division process described above, constituting the basic building blocks of a STIT tessellation as discussed in detail below. We shall provide explicit as well as asymptotic variance expressions for these parameters of the random tessellation $Y(t,W_R)$ as $R$ tends to infinity for a sequence $W_R=RW$ of expanding convex windows. Further, we will find the covariance measures for random lower-dimensional face measures generated by $Y(t,W).$ Finally, we shall also give the corresponding convergence in law statements, obtaining non-Gaussian limits for the studied case $d > 2$ as opposed to the classical Gaussian limits arising for $d=2$, see \cite{ST2}.\\ \\ The paper is organized as follows: In Section \ref{secSTIT} we define STIT tessellations,
 specializing to the stationary and isotropic set-up in the focus of this paper, and we discuss their basic
 properties that are needed in our arguments. Next, in Section \ref{sec2ndORD} we calculate the variances of the cumulative intrinsic volumes of all orders for STIT tessellations. This includes both, exact formulae and asymptotic analysis upon letting the window size grow to infinity. Further, in Section \ref{sec2ndMEAS} we extend the second-order analysis to the level of lower-dimensional face measures induced by STIT tessellations, thus taking into account not only the numeric characteristics but also the spatial profile of the STIT face processes. Finally, in Section \ref{secCLT} we develop the corresponding functional limit theory with non-Gaussian limit processes for $d > 2.$ In order to keep the paper self-contained, we will recall important faces from \cite{ST} and sketch some of their proofs for the readers convenience.

\section{STIT Tessellations}\label{secSTIT}

 The purpose of this Section is to provide a short self-contained discussion of STIT tessellations
 in ${\Bbb R}^d$ as studied in this paper and to summarize their basic properties for easy reference. We will restrict to the stationary and isotropic case in Subsection
 \ref{subsecMNWC} below, specializing to the scope of the paper. The general reference
 for these -- by now classical -- properties  throughout this section is \cite{NW05}. Next, in Subsection \ref{subsecMART}
 we will discuss the martingale tools developed in \cite{ST} and underlying our present theory. Finally,
 in Subsection \ref{subsec1stORD}
 we provide certain useful mean value relationships for intrinsic volumes in context of STIT tessellations.

\subsection{The MNW-Construction and Basic Properties}\label{subsecMNWC}

\begin{figure}[t]
 \begin{center}
  \includegraphics[width=7cm]{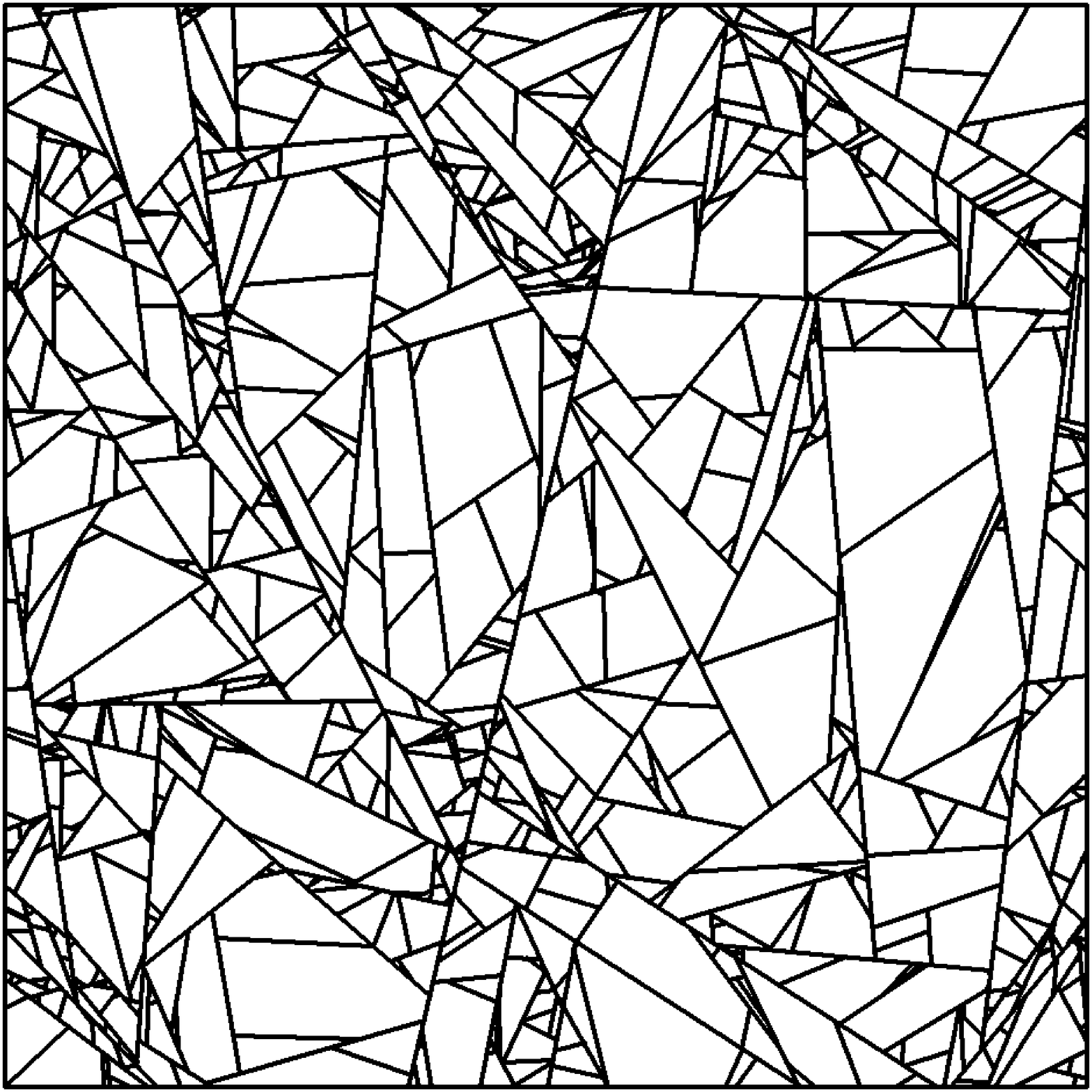}
  \includegraphics[width=7cm]{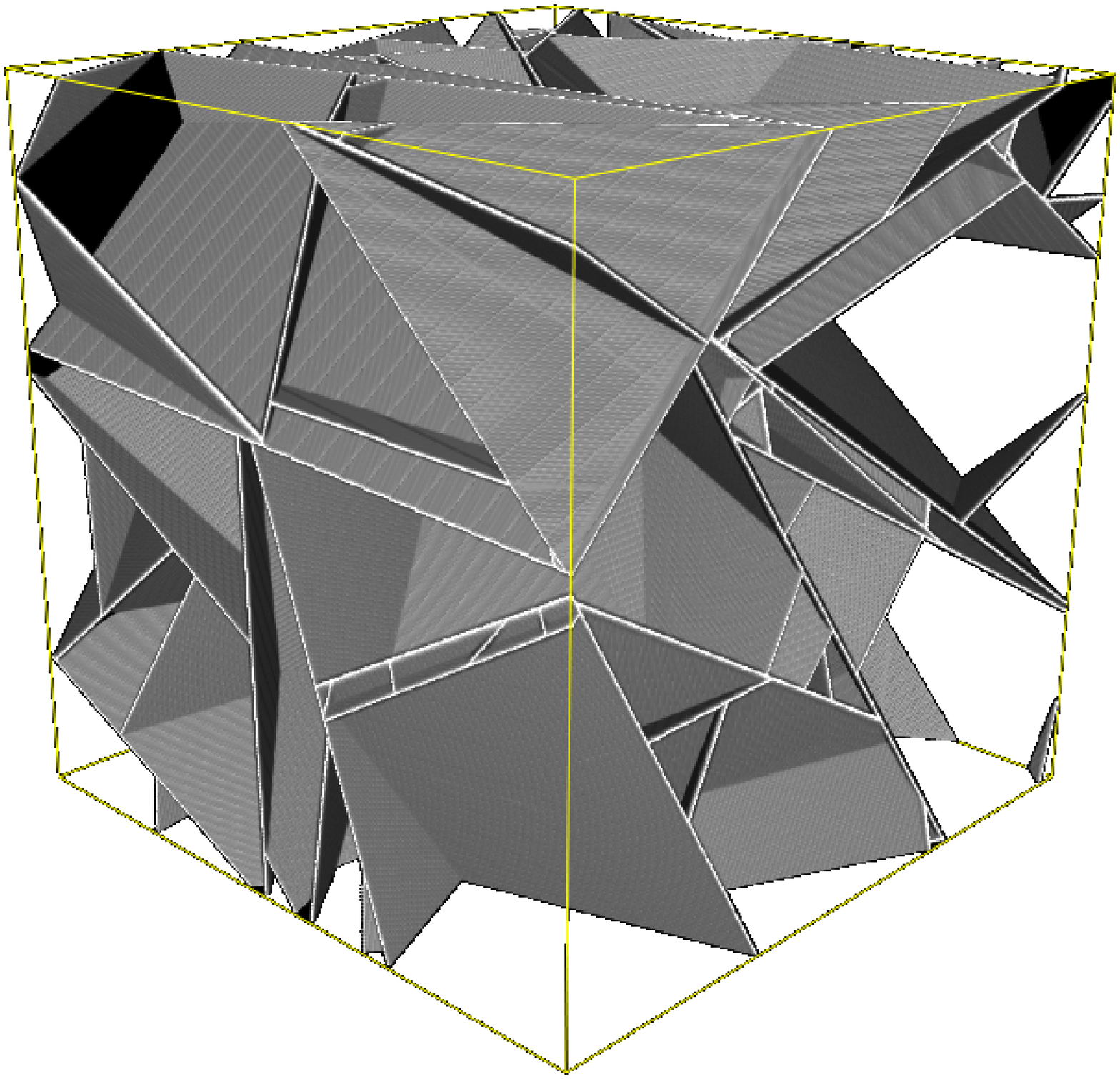}
  \caption{Realizations of a planar and a spatial stationary and isotropic STIT tessellation}
  \label{Figure}
 \end{center}
\end{figure}

We start with a compact and convex polytope $W\subset{\Bbb R}^d$ in which our construction is carried out and denote by $\Lambda$ the standard isometry-invariant measure on the space $\cal H$ of (affine) hyperplanes in
${\Bbb R}^d$ normalized so as to induce unit surface intensity on ${\cal H}$ (note that in our earlier papers \cite{ST} and \cite{ST2} this measure has been denoted by $\Lambda_{iso}$). We call $\Lambda$ the \textit{driving measure} of the construction. Assign now to $W$ an exponentially random lifetime with parameter $\Lambda([W])$, where $$[W]:=\{H\in{\cal H}:\ H\cap W\neq\emptyset\}$$ is the set of hyperplanes hitting $W$. Upon expiry of this random life time, a random hyperplane is chosen according to the distribution $\Lambda([W])^{-1}\Lambda(\cdot\cap[W])$, is introduced in $W$ and is chopped off by its boundary. This is, the window $W$ splits into the two polyhedral sub-cells $W^+$ and $W^-$ that are separated by the introduced hyperplane piece. The construction continues now recursively and independently in $W^+$ and $W^-$ and is stopped if some previously fixed deterministic time threshold $t>0$ is reached. Our assumptions ensure that the cells of the tessellation constructed until time $t>0$ are convex polyhedra in $W$ with probability one. They are denoted by $\Cells(Y(t,W))$ and we denote by $Y(t,W)$ the random closed set in $W$ that consists of the union of cell-boundaries of cells constructed until time $t$, see Figure \ref{Figure}. The construction of $Y(t,W)$ is referred to as the \textit{MNW-construction} after the names of its inventors \--- Mecke, Nagel and Weiss \--- and the tessellation $Y(t,W)$ itself is called a random \textit{STIT tessellation}. The abbreviation STIT comes from the crucial property enjoyed by $Y(t,W),$ namely that of being \underline{st}able under \underline{it}eration, for which we refer to \cite{MNW} or \cite{NW05}. We call the cell separating $(d-1)$-dimensional faces the \textit{$(d-1)$-dimensional maximal polytopes} and denote the collection of all such polytopes of $Y(t,W)$ by $\MaxFacets_{d-1}(Y(t,W))$. Moreover, we introduce the set $\MaxFaces_j(Y(t,W))$ of $j$-dimensional faces of $(d-1)$-dimensional maximal polytopes of $Y(t,W)$ by $$\MaxFaces_j(Y(t,W))=\bigcup_{f\in\MaxFacets_{d-1}(Y(t,W))}\Faces_j(f),$$ where by $\Faces_j(f)$ we mean the set of all $j$-dimensional faces of the $(d-1)$-dimensional polytope $f$, $0\leq j\leq d-1$. The maximal polytopes are often referred to as \textit{I-polytopes} \--- this terminology
originates from the historically first considered particular case $d=2$ where the maximal polytopes (which are just line segments in this case) assumed shapes
similar to the literal $I.$\\
It is an important observation that the spatio-temporal construction of $Y(t,W)$ satisfies the Markov property in time parameter $t$, which means that $$Y(t+s,W)=Y(t,W)\boxplus Y(s,W),$$ where $\boxplus$ denotes the operation of iteration of tessellations, see \cite{NW05}, \cite{MNW} or \cite{ST}. To make the paper more self-contained, we recall now the definition of the meaning of the operation $\boxplus$. To this end, we regard the tessellation $Y(t,W)$ as a frame or primary tessellation and associate with each cell $c\in\Cells(Y(t,W))$ an i.i.d. copy $Y_c(s,W)$ of $Y(s,W)$. Another tessellation $\tilde{Y}(t,s,W)$ of $W$ is now defined as $$\tilde{Y}(t,s,W):=Y(t,W)\cup\bigcup_{c\in\Cells(Y(t,W))}(c\cap Y_c(s,W)),$$ i.e. we consider the union of $Y(t,W)$ and the cut-outs of $Y_c(s,W)$'s within the cells $c$ of the primary tessellation (note that in the definition we have used the interpretation of a tessellation as a random closed subset of $W$). We say that $\tilde{Y}(t,s,W)$ is the \textit{iteration} of $Y(t,W)$ with $Y(s,W)$. The remarkable property of our tessellations constructed by the MNW-process is that the outcome $\tilde{Y}(s,t,W)$ coincides in law with $Y(t+s,W)$, i.e. with the continuation of the MNW-construction until time $t+s$. Thus, it is the same either to continue the MNW-construction from $t$ until time $t+s$ or to perform at time $t$ an iteration of $Y(t,W)$ with $Y(s,W)$.\\ The local properties established so far, can be extended to the whole space, since the random tessellations $Y(t,W)$ are \textit{consistent} in $W$, by which we mean that for any fixed $t>0$ and $W\subset W'\subset{\Bbb R}^d$ it holds $$Y(t,W')\cap W \overset{D}{=} Y(t,W),$$ where $\overset{D}{=}$ stands for equality in distribution. This implies -- in view of the consistency theorem \cite[Thm. 2.3.1]{SW} -- the existence of the whole-space tessellation $Y(t)$ such that $Y(t,W) \overset{D}{=} Y(t) \cap W$ for each compact convex $W\subset{\Bbb R}^d$.\\ It directly follows from the Markov property that the random tessellations $Y(t)$ are stable with respect to the operation $\boxplus$, i.e. $$Y(t)\overset{D}{=}n(\underbrace{Y(t)\boxplus\ldots\boxplus Y(t)}_{n}),\ \ \ \ \ n\in{\Bbb N},$$ where $n(\cdot)$ means the dilation with a factor $n$. This property explains also the abbreviation STIT, because the last equation is a classical probabilistic stability relation. The random tessellations $Y(t)$ share another important property, namely that the \textit{intersection} of $Y(t)$ with a $j$-dimensional plane $E_j$, $1\leq j\leq d-1$, is again an iteration stable random tessellation. More precisely, we have 
\begin{equation}\label{STITSECTION}
Y(t)\cap E_j \overset{D}{=} Y(\gamma_jt,E_j)
\end{equation}
and the integral-geometric constant $\gamma_j$ is given by 
\begin{equation}\label{GAMMAJ}
 \gamma_j = {\Gamma\left({j+1\over 2}\right)\Gamma\left({d\over 2}\right)\over\Gamma\left({j\over 2}\right)\Gamma\left({d+1\over 2}\right)},
\end{equation}
see Eq. (3.29T) in \cite{Miles}. Moreover, it is easy to see from the properties of the capacity functional of $Y(t)$, compare with \cite[Lem. 5(ii)]{NW05}, that STIT tessellations have the following \textit{scaling property}: \begin{equation} tY(t)\overset{D}{=}Y(1),\label{STITSCALING}\end{equation} i.e. the tessellation $Y(t)$ has, upon rescaling with the factor $t,$ the same distribution as $Y(1)$, that is the STIT tessellation with surface intensity $1$.\\ We close this section by mentioning that the random tessellations $Y(t)$ have another important property, namely that $Y(t)$ has \textit{Poisson typical cells} or Poisson cells for short, see Lemma 3 in \cite{NW05} or part (1) of Proposition \ref{PropXX} below. This is to say, the typical cell of $Y(t)$ has the same distribution as the typical cell of a stationary and isotropic Poisson hyperplane tessellation with surface density $t$ (see Chap. 10.3 in \cite{SW} for references about the classical Poisson hyperplane model). In particular, this fact combined with the intersection property from above shows that the typical cell of the lower dimensional STIT tessellation $Y(t)\cap E_j$, with $E_j$ as above, has the same distribution as the typical cell of a Poisson hyperplane tessellation in ${\Bbb R}^j$ with surface density $\gamma_jt$. More formally, we define the cell intensity measure ${\Bbb M}^{Y(t,W)}$ of the STIT tessellation $Y(t,W)$ by 
\begin{equation}\label{CellIntensityMeasure}{\Bbb M}^{Y(t,W)} := {\Bbb E}\sum_{c \in \Cells(Y(t,W))} \delta_c\end{equation}
and its $(d-1)$-dimensional maximal polytope intensity measure ${\Bbb F}_{d-1}^{Y(t,W)}$ by 
\begin{equation}\label{FaceIntensityMeasure}{\Bbb F}_{d-1}^{Y(t,W)} = {\Bbb E}\sum_{f \in \MaxFacets_{d-1}(Y(t,W))}\delta_f,\end{equation} where $\delta_{(\cdot)}$ stands for the uni mass Dirac measure concentrated at $(\cdot)$. Moreover, we let ${\Bbb M}^{\PHT(t,W)}$ be the cell intensity measure and ${\Bbb F}_{d-1}^{\PHT(t,W)}$ be the $(d-1)$-face intensity measure of a Poisson hyperplane tessellation $\PHT(t,W)$ within $W\subset{\Bbb R}^d$ having intensity measure $t\Lambda$, which are defined similarly to ${\Bbb M}^{Y(t,W)}$ and ${\Bbb F}_{d-1}^{Y(t,W)}$ above (clearly, the $(d-1)$-dimensional maximal polytopes have to be replaced by the set of $(d-1)$-faces of $\PHT(t,W)$). These definitions bring us in the position to reformulate special cases of Theorems 1 and 2 from \cite{ST} adapted to our later purposes:
\begin{proposition}\label{PropXX} We have
$$\text{(a)}\ \ {\Bbb M}^{Y(t,W)}={\Bbb M}^{\PHT(t,W)}\ \ \ \ \text{and}\ \ \ \ \text{(b)} \ \ {\Bbb F}_{d-1}^{Y(t,W)}=\int_0^t{1\over s}{\Bbb F}_{d-1}^{\PHT(s,W)}ds.$$

\end{proposition}

\subsection{Martingales in the MNW-Construction}\label{subsecMART}

As already noted above, the MNW-construction of iteration stable random tessellations $Y(t,W)$ in finite volumes $W\subset{\Bbb R}^d$ enjoys a Markov property in the continuous
time parameter $t.$ In our previous work \cite{ST} we have used this fact combined with the classical theory of martingale problems for pure jump Markov processes to construct a class of natural martingales associated
to the MNW-process. In this paper we only need a part of that theory. To formulate it we let $\phi$ be a
bounded and measurable functional on the space of $(d-1)$-polytopes in ${\Bbb R}^d$ and we denote by $\MaxFacets_{d-1}(Y)$ the collection of $(d-1)$-dimensional maximal polytopes of a tessellation $Y$ standing for a generic realization of
$Y(t,W)$ for some $t > 0.$ Write
\begin{equation}\label{SIGMADEF}
 \Sigma_\phi=\Sigma_\phi(Y):=\sum_{f\in\MaxFacets_{d-1}(Y)}\phi(f).
\end{equation}
Any hyperplane $H\in[W]$ hitting the window $W$ is tessellated by the intersection with $Y$ and we denote by $\Cells(Y\cap H)$
the set of all $(d-1)$-dimensional cells of this tessellation and introduce 
\begin{equation}\label{ADEF}
 A_\phi(Y):=\int_{[W]}\sum_{f\in\Cells(Y\cap H)}\phi(f)\Lambda(dH).
\end{equation}
It is also convenient to introduce the bar notation for centered versions of these quantities with
$Y = Y(t,W),$ that is to say $\bar{\Sigma}_{\phi}(Y(t,W)) := \Sigma_{\phi}(Y(t,W)) -
 {\Bbb E}\Sigma_{\phi}(Y(t,W))$
and likewise $\bar{A}_{\phi}(Y(t,W)) := A_{\phi}(Y(t,W)) - {\Bbb E}A_{\phi}(Y(t,W)).$
With this notation, in view of the results developed in \cite{ST}, we have
\begin{proposition}\label{MartProp}
 For bounded and measurable functionals $\phi$ and $\psi$ on the space 
of $(d-1)$-polytopes in ${\Bbb R}^d$, the stochastic processes
\begin{equation}\label{MART1}
\Sigma_\phi(Y(t,W))-\int_0^tA_\phi(Y(s,W))ds
\end{equation}
and
\begin{eqnarray}
\nonumber & & \bar\Sigma_\phi(Y(t,W))\bar\Sigma_\psi(Y(t,W))-\int_0^tA_{\phi\psi}(Y(s,W))ds\\
 &-&  \int_0^t[\bar A_\phi(Y(s,W))\bar\Sigma_\psi(Y(s,W))+\bar A_\psi(Y(s,W))\bar\Sigma_\phi(Y(s,W))]ds\label{MART2}
\end{eqnarray}
are martingales with respect to the filtration $\Im_t$ induced by $(Y(s,W))_{0\leq s\leq t}$.
\end{proposition}
\paragraph{Sketch of a Proof:} In order to make the paper self-contained, we give here the main idea of the proof of Proposition \ref{MartProp}, although there is some overlap with \cite{ST}.\\ At first, it is a direct consequence of the MNW-construction that the generator ${\Bbb L}$ of the pure-jump Markov process $Y(t,W)$, $t>0$, is given by $${\Bbb L}F(Y)=\int_{[W]}\sum_{f\in\Cells(Y\cap H)}[F(Y\cup f)-F(Y)]\Lambda(dH),$$ where $Y$ stands for some instant of $Y(t,W)$ and $F$ is a bounded and measurable function on the space of tessellations of $W$. Applying now the classical Dynkin formula (see for example Lemma 19.21 in \cite{Kallenberg}) with $F=\Sigma_\phi(Y)$ gives the martingale property of the random process (\ref{MART1}). (It should be mentioned that so defined $F$ need not to be bounded in general. However, this technical difficulty can be overcome with a suitable localization argument.)\\ For the second statement we consider the time-augmented random process $(Y(t,W),t)$, $t>0$, which has generator ${\Bbb L}'$ given by $${\Bbb L}'G(Y,t)=[{\Bbb L}G(\cdot,t)](Y)+\left[{\partial\over\partial t}G(Y,\cdot)\right](t)$$ for appropriate functions $G(Y,t)$. Now, using Dynkin's formula now for the product process $(Y(t,W),t)$ and with $G(Y,t):=(\Sigma_\phi(Y)-{\Bbb E}\Sigma_\phi(Y))^2=\bar\Sigma_\phi^2(Y)$, $Y=Y(t,W)$, gives (again after a suitable localization argument) the martingale property of the random process $$\bar\Sigma_\phi^2(Y(t,W))-\int_0^tA_{\phi^2}(Y(s,W))ds-2\int_0^t\bar A_{\phi}(Y(s,W))\bar\Sigma_\phi(Y(s,W))ds.$$ For another functional $\psi$ we can apply the latter property for $\phi+\psi$ and $\phi-\psi$ and subtract the two results, which shows the martingale property of the random process (\ref{MART2}).\hfill $\Box$

\subsection{Mean Values for Intrinsic Volumes}\label{subsec1stORD}

 In this subsection we discuss certain basic first-order consequences of Proposition \ref{MartProp} to
 be of use for our future reference. We will write $V_j$ for the intrinsic volume of order $j$ with $j=0,\ldots,d-1$ (for the standard definition of these functionals we refer to \cite{SW} and the references cited therein). Further, write $$F_j(Y):=\sum_{c \in \Cells(Y(t,W))}V_j(c).$$ Note now that whenever a new facet $f$ splits a cell $c$ into $c^+$ and $c^-$ of $Y(t,W)$ giving rise to a new tessellation $Y'$, we have $$F_j(Y') - F_j(Y(t,W)) = V_j(c^+)+V_j(c^-)-V_j(c) = V_j(c^+\cap c^-) = V_j(f),$$ since $V_j$ has the valuation property. Consequently, constructing the tessellation $Y=Y(t,W)$ by successive cell splits we easily get
\begin{equation}\label{FjSj}
  F_j(Y) = \Sigma_{V_j}(Y) + V_j(W).
\end{equation}
It is our aim to relate $F_j(Y)$ and $F_{j+1}(Y)$. For this purpose we use a special case of \textit{Crofton's formula} from classical integral geometry, which reads \begin{equation}\int_{[K]}V_j(K\cap H)\Lambda(dH)=\gamma_{j+1}V_{j+1}(K),\label{EQCROFTON}\end{equation} where $K\subset{\Bbb R}^d$ is a convex body and where the constant $\gamma_j$ is given by (\ref{GAMMAJ}), see \cite[Thm. 5.1.1]{SW}. Applying this formula to $F_j(Y)$ yields 
\begin{eqnarray}
\nonumber A_{V_j}(Y) &=& \int_{[W]}\sum_{f\in\Cells(Y\cap H)}V_j(f)\Lambda(dH)=\int_{[W]}\sum_{c\in\Cells(Y)}V_j(c\cap H)\Lambda(dH)\\
 &=& \sum_{c\in\Cells(Y)}\gamma_{j+1}V_{j+1}(c)=\gamma_{j+1} F_{j+1}(Y).\label{AVJTECHNICALTOOL}
\end{eqnarray}
To streamline the notation, we
shall write  $\Sigma_{V_j;s} := \Sigma_{V_j}(Y(s,W))$ and $A_{V_j;s} := A_{V_j}(Y(s,W))$ below. 
Note that, upon combining (\ref{FjSj}) and (\ref{AVJTECHNICALTOOL}), it follows from (\ref{MART1})
with $\phi = V_j$ that
\begin{equation}\label{MartExpl1ord}
 \Sigma_{V_j;t} - \gamma_{j+1} \int_0^t \Sigma_{V_{j+1};s} ds - t \gamma_{j+1} V_{j+1}(W) 
\end{equation} 
are all $\Im_t$-martingales for $j=0,\ldots,d-1.$  In particular, bearing in mind that $\Sigma_{V_d;t} \equiv 0$
we conclude that $\Sigma_{V_{d-1};t} - t \gamma_d V_d(W)$ is a $\Im_t$-martingale, in particular ${\Bbb E}\Sigma_{V_{d-1};t}=tV_d(W)$, where we have used $\gamma_d=1$. The latter equation is extended by 
\begin{proposition}\label{propoFIRSTORDER} For $j\in\{0,1,\ldots,d-1\}$ we have $${\Bbb E}\Sigma_{V_j;t}=\left(\prod_{i=1}^{d-j-1}\gamma_{j+i}\right){t^{d-j}\over(d-j)!}V_d(W)+\sum_{i=1}^{d-j-1}\left(\prod_{k=1}^{i}\gamma_{j+k}\right){t^{i}\over i!}V_{j+i}(W)$$ with the convention that $\prod_{i=1}^0\ldots\equiv 1$ and $\sum_{i=1}^0\ldots\equiv 0$.
\end{proposition}
\paragraph{Proof:} 
Taking expectations in (\ref{MartExpl1ord}) we get
$${\Bbb E}\Sigma_{V_j;t}=\int_0^t\gamma_{j+1}{\Bbb E}\Sigma_{V_{j+1};s_1}ds_1+\gamma_{j+1}tV_{j+1}(W).$$ Continuing recursively by applying (\ref{FjSj}) and (\ref{AVJTECHNICALTOOL}) we end up with
\begin{eqnarray}
\nonumber {\Bbb E}\Sigma_{V_j;t} &=& \int_0^t\cdots\int_0^{s_{d-j-2}}\gamma_{j+1}\cdots\gamma_{d-1}{\Bbb E}\Sigma_{V_{d-1};s_{d-j-1}} ds_{d-j-1}\cdots ds_1\\
\nonumber & & +\gamma_{j+1}\cdots\gamma_{d-1}{t^{d-j-1}\over(d-j-1)!}V_{d-1}(W)+\ldots+\gamma_{j+1}tV_{j+1}(W).
\end{eqnarray}
However, the definitions of $Y(s,W)$ and $V_{d-1}$ imply that ${\Bbb E}\Sigma_{V_{d-1};s}=sV_d(W)$ for any $0 \leq s \leq t$. Thus,
\begin{eqnarray}
\nonumber & & \int_0^t\cdots\int_0^{s_{d-j-2}}\gamma_{j+1}\cdots\gamma_{d-1}{\Bbb E}\Sigma_{V_{d-1};s_{d-j-1}}ds_{d-j-1}\cdots ds_1\\
\nonumber &=& \gamma_{j+1}\cdots\gamma_{d-1}V_d(W)\int_0^t\cdots\int_0^{s_{d-j-2}}s_{d-j-1}ds_{d-j-1}\cdots ds_1\\
\nonumber &=& \left(\prod_{i=1}^{d-j-1}\gamma_{j+i}\right){t^{d-j}\over(d-j)!}V_d(W)
\end{eqnarray}
and the result follows immediately.\hfill $\Box$\\ \\
It is interesting to compare the mean value formula from the last proposition with that from \cite{ST}, Section 3.2. We denote by $\overline{\varphi}_j(Y(t))$ the density of the $j$-th intrinsic volume of the collection of $(d-1)$-dimensional maximal polytopes of $Y(t)$, i.e. 
\begin{eqnarray}
\nonumber \overline{\varphi}_j(Y(t)) &=& \lim_{r\rightarrow\infty}{1\over r^dV_d(W)}{\Bbb E}\sum_{f\in\MaxFacets_{d-1}(Y(t,rW))}V_j(f)\\
\nonumber &=& \lim_{r\rightarrow\infty}{1\over r^dV_d(W)}{\Bbb E}\Sigma_{V_j}(Y(t,rW))
\end{eqnarray}
for arbitrary compact windows $W\subset{\Bbb R}^d$ with $V_d(W)>0$, see \cite[Chap. 4.1]{SW}.
 Using now Proposition \ref{propoFIRSTORDER} we get $$\overline{\varphi}_j(Y(t))=\left(\prod_{i=1}^{d-1-j}\gamma_{i+j}\right){t^{d-j}\over(d-j)!},$$ because of the homogeneity of intrinsic volumes. On the other hand, we have shown in \cite{ST} that \begin{equation}\overline{\varphi}_j(Y(t))={d\choose j}\left({\kappa_{d-1}\over d\kappa_d}\right)^{d-j}{\kappa_d\over\kappa_j}t^{d-j},\label{XXXEQMEAN}\end{equation} where here and in the sequel we denote by $\kappa_i$, $0\leq i\leq d$, the volume of the $i$-dimensional unit ball. Indeed, these two values are identical, since
\begin{eqnarray}
\nonumber {d\choose j}\left({\kappa_{d-1}\over d\kappa_d}\right)^{d-j}{\kappa_d\over\kappa_j} &=& {\Gamma\left({d+1\over 2}\right)\over\Gamma\left({j+1\over 2}\right)}\left({\Gamma\left({d\over 2}\right)\over\Gamma\left({d+1\over 2}\right)}\right)^{d-j}{1\over(d-j)!}\\
\nonumber &=& \left({\Gamma\left({d\over 2}\right)\over\Gamma\left({d+1\over 2}\right)}\right)^{d-1-j}{\Gamma\left({d\over 2}\right)\over\Gamma\left({j+1\over 2}\right)}{1\over(d-j)!}\\
\nonumber &=& \left(\prod_{i=1}^{d-1-j}\gamma_{i+j}\right){1\over(d-j)!}.
\end{eqnarray}
\begin{remark}\label{Rema1} Our first-order formula from the last proposition contains the extra boundary correction term $\displaystyle\sum_{i=1}^{d-j-1}\left(\prod_{k=1}^{i}\gamma_{j+k}\right){t^{i}\over i!}V_{j+i}(W)$ in contrast to our mean value formula (\ref{XXXEQMEAN}) from \cite{ST}. These additional terms come from the fact that we consider maximal polytopes possibly chopped off by the boundary of the domain $W$ rather than the points of an
associated center function for full facets in ${\Bbb R}^d.$ Thus, it may happen for instance that in two neighbouring regions one observes two facets which can coalesce when putting these regions together into one volume, whence
the lower-order finite volume corrections arise.
\end{remark}

\section{Second Order Characteristics for Intrinsic Volumes}\label{sec2ndORD}
 In this Section we develop a full second-order theory for intrinsic volumes of stationary and isotropic random
 STIT tessellations. This is first done exactly in Subsection \ref{subsecExact} for $Y(t,W),\; t>0$ inside
 a bounded convex observation window $W\subset{\Bbb R}^d$. Then, in Subsection \ref{subsecAsympt}
 we derive the corresponding asymptotic expressions for $Y(t,W_R),\; W_R := RW,\; R \to \infty.$

\subsection{Exact Expressions}\label{subsecExact}
 To proceed with second-order calculations we shall use the notation already introduced in Subsection \ref{subsec1stORD} above. Observe first that in view of (\ref{FjSj})
 the relation (\ref{AVJTECHNICALTOOL}) simplifies and we obtain
\begin{equation}\label{BarAj}
  \bar{A}_{V_j}(Y) = \gamma_{j+1} \bar{\Sigma}_{V_{j+1}}(Y).
\end{equation}
Recalling that $\Sigma_{V_j;s} := \Sigma_{V_j}(Y(s,W))$ and $A_{V_j;s} := A_{V_j}(Y(s,W)),$
putting likewise $A_{V_i V_j;s} := A_{V_i V_j}(Y(s,W))$ 
and using (\ref{MART2}) in Proposition \ref{MartProp} with $\phi = V_i$ and $\psi = V_j$, we see that
 $$ \bar{\Sigma}_{V_i;t} \bar{\Sigma}_{V_j;t} - \int_0^t A_{V_i V_j;s} ds
  - \int_0^t [\bar{A}_{V_i;s} \bar{\Sigma}_{V_j;s} + \bar{A}_{V_j;s} \bar{\Sigma}_{V_i;s}] ds $$
 is a martingale with respect to the filtration $\Im_t$ induced by $(Y(t,W))_{0 \leq s \leq t}$.
 Substituting (\ref{BarAj}), this implies that
 \begin{equation}\label{PrincipalMart}
 \nonumber \bar{\Sigma}_{V_i;t} \bar{\Sigma}_{V_j;t} - \int_0^t A_{V_i V_j;s} ds
  - \int_0^t [\gamma_{i+1} \bar{\Sigma}_{V_{i+1};s} \bar{\Sigma}_{V_j;s}
  + \gamma_{j+1} \bar{\Sigma}_{V_i;s} \bar{\Sigma}_{V_{j+1};s}] ds
 \end{equation}
 is an $\Im_t$-martingale as well. This is a crucial formula because, upon taking 
 expectations, it allows us to express $\Cov(\Sigma_{V_i;t},\Sigma_{V_j;t})$ 
 in terms of corresponding covariances with indices $(i,j+1)$ and $(j,i+1).$ In other words,
 we get
 \begin{equation}\label{PrincipalCov}
  \Cov(\Sigma_{V_i;t},\Sigma_{V_j;t}) = \int_0^t {\Bbb E} A_{V_i V_j;s} ds 
  + \int_0^t [\gamma_{i+1} \Cov(\Sigma_{V_{i+1};s},\Sigma_{V_j;s}) +
     \gamma_{j+1} \Cov(\Sigma_{V_i;s},\Sigma_{V_{j+1};s})] ds.
 \end{equation}
It is important to observe that this recursion terminates because $\Sigma_{V_d;t} = 0$, which allows us to provide an explicit expression for the covariances $\Cov(\Sigma_{V_i;t},\Sigma_{V_j;t})$. To obtain the desired formula, denote
 \begin{equation}\label{ItInt}
  {\cal I}^n(f;t) := \int_0^t \int_0^{s_1} \ldots \int_0^{s_{n-1}} f(s_n) ds_n d s_{n-1} \ldots ds_1 
    = \frac{1}{(n-1)!} \int_0^t (t-s)^{n-1} f(s) ds
 \end{equation}
 for any $f : [0,t] \to {\Bbb R}$ and $n\in{\Bbb N}$ for which the iterated integral exists. Then we claim that for $k,l\in\{0,\ldots,d-1\}$ it holds
 $$
  \Cov(\Sigma_{V_{d-1-k};t},\Sigma_{V_{d-1-l};t}) = {{k+l}\choose k} \left( \prod_{i=1}^k \gamma_{d-i} \right)
  \left( \prod_{j=1}^l \gamma_{d-j} \right) {\cal I}^{k+l}(\Var(\Sigma_{V_{d-1};(\cdot)});t) +
 $$
 \begin{equation}\label{AuxEq}
   \sum_{\stackrel{0\leq m \leq k,\; 0 \leq n \leq l}{(m,n) \neq (0,0)}}
   {{k+l-m-n} \choose k-m}  
     \left( \prod_{i=m+1}^k \gamma_{d-i} \right) \left( \prod_{j=n+1}^l \gamma_{d-j} \right) 
     {\cal I}^{k+l-m-n+1}({\Bbb E}A_{V_{d-1-m} V_{d-1-n};(\cdot)};t) 
 \end{equation}
 with the convention that ${0\choose 0} = 1$ and $\prod_{k+1}^k \ldots = \prod_{l+1}^l \ldots \equiv 1.$
 Whereas this can be readily verified by a straightforward induction in view of (\ref{PrincipalCov}),
 there is a more natural way to see it. In fact, in the course of recursive applications of
 (\ref{PrincipalCov}) the covariance $c_{k,l} := \Cov(\Sigma_{V_{d-1-k}},\Sigma_{V_{d-1-l}})$
 is represented in terms of $c_{k-1,l}$ and $c_{k,l-1}$ which can be interpreted in terms
 of lattice walks on pairs of indices from $(k,l)$ to $(0,0)$ where only steps $(i,j) \to (i-1,j)$ and
 $(i,j) \to (i,j-1)$ are allowed, each step $(i,j) \to (i-1,j)$ involving multiplication by $\gamma_{d-i}$
 plus integral iteration and each step $(i,j) \to (i,j-1)$ resulting in multiplication by $\gamma_{d-j}$
 plus integral iteration. There are $k+l\choose k$ such walks, whence the first line in (\ref{AuxEq})
 follows. The second line of (\ref{AuxEq}) is due to additional 
 $\int_0^t {\Bbb E}A_{V_{d-1-k+m'} V_{d-1-l+n'};s} ds$ terms
 in (\ref{PrincipalCov}), which are born at all times $(m',n')$ of the lattice walk discussed above,
 with $m'$ standing for the number of $(i,j) \to (i-1,j)$ steps and $n'$ for the number of $(i,j) \to (i,j-1)$
 steps. Note that no additional term is born at $(m'=k,n'=l)$ though, as it corresponds to having reached
 the (co)variance $\Cov(\Sigma_{V_{d-1}},\Sigma_{V_{d-1}})=\Var(\Sigma_{V_{d-1}})$, which does not get expanded any further.
 In (\ref{AuxEq}) we have substituted $m = k - m'$ and $n = l - n'$ for notational convenience.\\
 Applying now (\ref{PrincipalCov}) 
 for $i=j=d-1$ and using that $\Sigma_{V_d;t} = 0$ we see that 
 \begin{equation}\label{VarD1expr}
  \Var(\Sigma_{V_{d-1};t}) = {\cal I}^{1}({\Bbb E}A_{V^2_{d-1};(\cdot)};t), 
 \end{equation}
 whence we finally get from (\ref{AuxEq}) the following exact variance expression:
\begin{theorem}\label{thmconvexact} The covariance between the intrinsic volume processes
 $\Sigma_{V_{d-1-k};t}$ and $\Sigma_{V_{d-1-l};t}$ for $k,l\in\{0,\ldots,d-1\}$ of a stationary and isotropic random STIT tessellation $Y(t,W)$ is given by 
 $$  \Cov(\Sigma_{V_{d-1-k};t},\Sigma_{V_{d-1-l};t}) = $$
 \begin{equation}\nonumber\label{ExplicitCov}
   \sum_{m=0}^k \sum_{n=0}^l
   {{k+l-m-n} \choose k-m} 
     \left( \prod_{i=m+1}^k \gamma_{d-i} \right) \left( \prod_{j=n+1}^l \gamma_{d-j} \right) 
     {\cal I}^{k+l-m-n+1}({\Bbb E}A_{V_{d-1-m} V_{d-1-n};(\cdot)};t) 
 \end{equation}
  with the usual convention that ${0\choose 0} = 1$ and $\prod_{k+1}^k \ldots = \prod_{l+1}^l \ldots \equiv 1.$
\end{theorem}
 Unfortunately, we are not able to make the covariance in the above theorem any more explicit. This is due to the presence
 of mixed moments ${\Bbb E}A_{V_{d-1-m} V_{d-1-n};t},\; m,n > 0,$ whose evaluation is technically
 related to the problem of providing a general formula for joint moments of lower-order intrinsic volumes
 of $(d-1)$-dimensional Poisson cells (arising as hyper-planar sections of $Y(t,W)$, which are in addition possibly
 chopped off by the boundary of the window $W$) which is not currently available in required generality up to the best of
 our knowledge. Fortunately though, the offending mixed moments turn out to be of negligible order
 in large window size asymptotics and we are able to provide fully explicit asymptotic formulae in 
 Subsection \ref{subsecAsympt} below.\\
We would like to point out that in the special case $d=2$ our Theorem \ref{thmconvexact} reduces to Theorem 1 and Corollary 1 in \cite{ST2} where we have studied variances and central limit theory for maximal edge (I-segment) lengths and vertex counts and where we could give fully explicit exact variance formulae exploiting the particular features of the planar setting. In Corollary 2 ibidem we also have provided asymptotic expressions for these variances for sequences of growing observation windows. For this reason, in our asymptotic considerations below we will restrict to the case $d>2$ where the general asymptotic expressions are essentially different from those arising in the exceptional planar case. This is due to the variance dichotomy established in \cite{ST} and mentioned in Section
\ref{secIntro} above.

\subsection{Asymptotic Expressions for $d>2$}\label{subsecAsympt}

Let $W\subset{\Bbb R}^d$ be a compact convex window and consider the sequence $W_R:=RW$. We will write from now on $A^{W_R}_{\phi;t}$ instead of $A_\phi(Y(t,W_R))$ and likewise $\Sigma_{\phi;t}^{W_R}$ for $\Sigma_\phi(Y(t,W_R))$ in order to emphasize the dependence on $R$. Our main interest in this section is to derive from Theorem \ref{thmconvexact} an asymptotic expression for the covariances $\Cov(\Sigma_{V_{d-1-k};t}^{W_R},\Sigma_{V_{d-1-l};t}^{W_R})$ as $R\rightarrow\infty$. We start with the following

\begin{proposition}\label{PropAux1}
 For $k,l,m,n \in \{0,\ldots,d-1\},\; m \leq k,\; n \leq l,$ with $t$ fixed we have
 $$  {\cal I}^{k+l-m-n+1}({\Bbb E}A^{W_R}_{V_{d-1-m} V_{d-1-n};(\cdot)};t) = 
\left\{ \begin{array}{ll}
 O(R^{2(d-1)-m-n}), & \mbox{ if } m+n \leq d-3,\\
 O(R^d \log R), & \mbox{ if } m+n = d - 2,\\
 O(R^d), & \mbox{ if } m+n \geq d-1. \end{array} \right. 
$$
\end{proposition}

\paragraph{Proof:}  Observe first that, in view of (\ref{ItInt}), the integral ${\cal I}^{k+l-m-n+1}$ only involves integration with respect to the variable $s$ and does not affect the order in $R$. Thus, it is sufficient to prove
the Proposition for $k=m$ and $l=n$ which we shall henceforth assume without loss of generality. We have
$$ {\cal I}^1({\Bbb E}A^{W_R}_{V_{d-1-m} V_{d-1-n};(\cdot)};t) = S^{(1)}_{R;t} + S^{(2)}_{R;t} $$
where
$$ S^{(1)}_{R;t} = \int_{1/R}^{t}  {\Bbb E}A^{W_R}_{V_{d-1-m} V_{d-1-n};s} ds $$
and
$$ S^{(2)}_{R;t} = \int_0^{1/R} {\Bbb E}A^{W_R}_{V_{d-1-m} V_{d-1-n};s} ds. $$
To provide a bound for $S^{(1)}_{R;t}$ we need some additional notation. Write
$\varsigma_{m,n;s}$ for the common value of
$$ {\Bbb E}V_{d-1-m}(\TypicalCell(Y(s)\cap H)) V_{d-1-n}(\TypicalCell(Y(s)\cap H))$$
with $H$ ranging through hyperplanes in ${\Bbb R}^d$ and where $\TypicalCell({Y(s) \cap H})$ the typical cell of the sectional STIT tessellation $Y(s) \cap H$.
 Using the scaling property (\ref{STITSCALING}) and the homogeneity of the intrinsic volumes we readily get
\begin{equation}\label{Bound1}
 \varsigma_{m,n;s} = s^{m+n-2(d-1)} \varsigma_{m,n;1}.
\end{equation}
Indeed,
\begin{eqnarray}
\nonumber \varsigma_{m,n;s} &=& {\Bbb E}V_{d-1-m}(\TypicalCell(Y(s)\cap H)) V_{d-1-n}(\TypicalCell(Y(s)\cap H))\\
\nonumber &=& {\Bbb E}V_{d-1-m}\left({1\over s}\TypicalCell(sY(s)\cap H)\right)V_{d-1-n}\left({1\over s}\TypicalCell(sY(s)\cap H)\right)\\
\nonumber &=& \left({1\over s}\right)^{d-1-m+d-1-n}{\Bbb E}V_{d-1-m}(\TypicalCell(Y(1)\cap H)) V_{d-1-n}(\TypicalCell(Y(1)\cap H))\\
\nonumber &=& s^{m+n-2(d-1)}\varsigma_{m,n;1}.
\end{eqnarray}
To proceed, write, recalling (\ref{ADEF}),
 \begin{eqnarray}
\nonumber {\Bbb E}A^{W_R}_{V_{d-1-m} V_{d-1-n};s} &=& {\Bbb E}\int_{[W_R]}\sum_{f\in\Cells(Y(s,W_R)\cap H)}V_{d-1-m}(f)V_{d-1-n}(f)\Lambda(dH)\\
\nonumber &\leq&
 {\Bbb E}\int_{[W_R]}\sum_{f\in\Cells(Y(s)\cap H),\; f \cap W_R \neq \emptyset}
  V_{d-1-m}(f)V_{d-1-n}(f)\Lambda(dH)\\
\nonumber &=& \int_{[W_R]} \frac{\varsigma_{m,n;s}}{\varsigma_{0,d-1;s}} 
 {\Bbb E}\Vol_{d-1}(\TypicalCell(Y(s) \cap H) \oplus (-(W_R \cap H))) \Lambda(dH),
\end{eqnarray}
where $\oplus$ stands for the usual Minkowski addition. Recalling that $\TypicalCell(Y(s) \cap H)$ is Poisson, see Section \ref{subsecMNWC}, we readily get for $s \geq 1/R$, $$   {\Bbb E}\Vol_{d-1}(\TypicalCell(Y(s) \cap H) \oplus (-(W_R \cap H))) = O(R^{d-1}) $$
and thus we conclude in view of (\ref{Bound1}) that, for $s \geq 1/R,$
$$ {\Bbb E}A^{W_R}_{V_{d-1-m} V_{d-1-n};s} = O(R \cdot s^{m+n-(d-1)} \cdot R^{d-1}) =
     O(R^d \cdot s^{m+n-(d-1)}). $$
Consequently,
\begin{equation}\label{BoundS1}
 S^{(1)}_{R;t} = O\left(R^d \int_{1/R}^t s^{m+n-(d-1)} ds\right) = \left\{ \begin{array}{ll}
 O(R^{2(d-1)-m-n}), & \mbox{ if } m+n \leq d-3,\\
 O(R^d \log R), & \mbox{ if } m+n = d - 2,\\
 O(R^d), & \mbox{ if } m+n \geq d-1. \end{array} \right. 
\end{equation} 
Next, we find a bound for $S^{(2)}_{R;t}.$ To this end, we note that during the time interval $[0,1/R]$
there are $O(1)$ cell splits within $W_R$ in the course of the MNW-construction and hence the expectations of sums 
$$\sum_{f \in \Cells(Y(s,W_R) \cap H)} V_{d-1-m}(f) V_{d-1-n}(f)$$ are of order $O(V_{d-1-m}(W_R)
 V_{d-1-n}(W_R)) = O(R^{2(d-1)-m-n})$ uniformly in hyperplanes $H \in [W_R]$ and in $s \in [0,1/R].$
Thus, using that $\Lambda([W_R]) = O(R)$ and recalling the definition (\ref{ADEF}) we obtain
\begin{equation}\label{BoundS2}
  S^{(2)}_{R;t} = O\left(R \cdot R^{2(d-1)-m-n} \cdot \int_0^{1/R} ds\right) = O(R^{2(d-1)-m-n}).
\end{equation}
Putting (\ref{BoundS1}) and (\ref{BoundS2}) together completes the proof of Proposition. \hfill $\Box$\\ \\
This implies that asymptotically, as $R\rightarrow\infty$, the terms appearing in Theorem \ref{thmconvexact} with $n,m>0$ are of order at most $O(R^{2d-3})$ and thus negligible compared with the leading $R^{d(d-1)}$-term. For this reason, the asymptotic behaviour of $\Cov(\Sigma_{V_{d-1-k};t},\Sigma_{V_{d-1-l};t})$ is dominated by the term with $m=n=0$, which is of order $\Theta(R^{2(d-1)}),$ where by $\Theta(\cdot)=O(\cdot)\cap\Omega(\cdot)$ we mean 
quantities bounded both from below and above by multiplicities of the argument in the usual Landau notation. For this dominating case Proposition \ref{PropAux1} is refined by
\begin{proposition}\label{PropAsympt}
 We have for $k,l\in\{0,\ldots,d-1\}$ $${\cal I}^{k+l+1}({\Bbb E}A^{W_R}_{V^2_{d-1};(\cdot)};t)={1\over(k+l)!}{d-1\over 2}t^{k+l}R^{2(d-1)}\int_W\int_W{1\over\left\|x-y\right\|^2}dxdy+O(R^{2d-3}).$$
\end{proposition}
\paragraph{Proof:} First, recall from Thm. 4 in \cite{ST} that the variance of the total surface area of $Y(t,W)$ equals
\begin{equation}\label{TOTALVAR}
\Var(\Sigma_{V_{d-1};t})={d-1\over 2}\int_W\int_W{1-e^{-{2\kappa_{d-1}\over d\kappa_d}t\left\|x-y\right\|}\over\left\|x-y\right\|^2}dxdy.
\end{equation}
The main argument for (\ref{TOTALVAR}) reads as follows: At first, we use the fact that STIT tessellations have Poisson typical cells and find
\begin{eqnarray}
\nonumber A_{V_{d-1}^2}(Y(s,W)) &=& \int_{[W]}\sum_{f\in\Cells(Y(s,W)\cap H)}V_{d-1}^2(f)\Lambda(dH)\\
\nonumber &=& \int_{[W]}\int_{W\cap H}\int_{W\cap H}{\bf 1}[x,y\ \text{are in the same cell of}\ Y(s,W)\cap H]dxdy\Lambda(dH)\\
\nonumber &=& \int_{[W]}\int_{W\cap H}\int_{W\cap H}e^{-s\Lambda([\overline{xy}])}dxdy\Lambda(dH),
\end{eqnarray}
where $\bf 1[\cdot]$ stands for the usual indicator function and $\overline{xy}$ for the line segment connecting $x$ and $y$. Using now (\ref{MART2}) with $\phi=\psi$, taking expectations and noting that the mixed term vanishes leads to
\begin{eqnarray}
\nonumber \Var(\Sigma_{V_{d-1};t}) = \int_0^tA_{V_{d-1}^2}(Y(s,W))ds &=& \int_0^t\int_{[W]}\int_{W\cap H}\int_{W\cap H}e^{-s\Lambda([\overline{xy}])}dxdy\Lambda(dH)\\
\nonumber &=& \int_{[W]}\int_{W\cap H}\int_{W\cap H}{1-e^{-t\Lambda([\overline{xy}])}\over\Lambda([\overline{xy}])}dxdy\Lambda(dH).
\end{eqnarray}
The latter integral can be transformed into (\ref{TOTALVAR}) by using an integral-geometric formula of Blaschke-Petkantschin-type and the Crofton formula (\ref{EQCROFTON}) -- see \cite{ST} for details.\\ Recall now (\ref{VarD1expr}) and write
\begin{equation}{\cal I}^{k+l+1}({\Bbb E}A^{W_R}_{V^2_{d-1};(\cdot)};t)={\cal I}^{k+l}(\Var(\Sigma_{V_{d-1};(\cdot)}^{W_R});t)=R^{2(d-1)}{\cal I}^{k+l}(\Var(\Sigma_{V_{d-1};(\cdot)R}^W);t),\label{eeee1}\end{equation} where the last equality follows by the scaling property (\ref{STITSCALING})
 of STIT tessellations. Thus, $${\cal I}^{k+l+1}({\Bbb E}A^{W_R}_{V^2_{d-1};(\cdot)};t)={R^{2(d-1)}\over(k+l-1)!}{d-1\over 2}\int_W\int_W\int_0^t(t-s)^{k+l-1}{1-e^{-{2\kappa_{d-1}\over d\kappa_d}Rs\left\|x-y\right\|}\over\left\|x-y\right\|^2}dsdxdy.$$ Now, the binomial theorem implies that there exists some constant $\tilde{c}>0$ such that
\begin{eqnarray}
\nonumber & & \int_0^t(t-s)^{k+l-1}{1-e^{-c\left\|x-y\right\|Rs}\over \left\|x-y\right\|^2}ds\\
\nonumber &=& {1\over(k+l)}{t^{k+l}c^{k+l}\left\|x-y\right\|^{k+l}R^{k+l}+\tilde{c}e^{-c\left\|x-y\right\|tR}+O(R^{k+l-1})\over \left\|x-y\right\|^{k+l+2}c^{k+l}R^{k+l}}\\
\nonumber &=& {t^{k+l}\over(k+l)}{1\over\left\|x-y\right\|^2}+O(1/R)
\end{eqnarray}
with $c={2\kappa_{d-1}\over d\kappa_d}$. In view of (\ref{eeee1}) this proves the desired result.\hfill $\Box$\\ \\
Consequently, by combining Propositions \ref{PropAux1} and \ref{PropAsympt}
 with Theorem \ref{thmconvexact} we get
\begin{corollary}\label{corasycov} Asymptotically, as $R\rightarrow\infty$ we have for $k,l\in\{0,\ldots,d-1\}$
 $$\Cov(\Sigma^{W_R}_{V_{d-1-k};t},\Sigma^{W_R}_{V_{d-1-l};t}) =$$
 $$ {k+l\choose k} \left( \prod_{i=1}^k \gamma_{d-i} \right)
      \left( \prod_{j=1}^l \gamma_{d-j} \right) 
      {1\over(k+l)!}{d-1\over 2}t^{k+l}R^{2(d-1)}E_2(W)+ O(R^{2d-3}),$$ where $E_2(W)$ denotes the $2$-energy of $W$, i.e. $$E_2(W)=\int_W\int_W{1\over\left\|x-y\right\|^2}dxdy.$$
\end{corollary}
The affine Blaschke-Petkantschin Formula \cite[Thm. 7.2.7]{SW} can be used to provide an integral-geometric expression for the measure-geometric energy functional $E_2$. In fact, we have $$E_2(W)={2\over(d-1)(d-2)}I_{d-1}(W),$$ where $I_{d-1}(W)$ denotes the $(d-1)$-st chord power integral of the convex body $W$ in the sense of \cite{SW}, page 363, i.e. $$I_{d-1}(W)={d\kappa_d\over 2}\int_{\cal L}\Vol_1^{d-1}(W\cap L)dL,$$ where ${\cal L}$ denotes the space of lines in ${\Bbb R}^d$ with invariant measure $dL$ and $\kappa_d$ is the volume of the $d$-dimensional unit ball. This means that in the asymptotic covariance formula from Corollary \ref{corasycov} the dependence on the geometry of $W$ is encoded by the non-additive $E_2(W)$ or equivalently by $I_{d-1}(W)$. We will from now on use the representation in terms of chord power integrals, as it allows an easier comparison with other tessellation models. Summarizing, this yields the following
\begin{corollary} Asymptotically, as $R\rightarrow\infty$, we have for $k,l\in\{0,\ldots,d-1\}$
$$\Cov(\Sigma^{W_R}_{V_{d-1-k};t},\Sigma^{W_R}_{V_{d-1-l};t})={1\over d-2}\left( \prod_{i=1}^k \gamma_{d-i} \right)
      \left( \prod_{j=1}^l \gamma_{d-j} \right){t^{k+l}\over k!l!}I_{d-1}(W)R^{2(d-1)}+O(R^{2d-3})$$ and for $k=l$ $$\Var(\Sigma_{V_{d-1-k};t}^{W_R})={1\over d-2}\left( \prod_{i=1}^k \gamma_{d-i} \right)^2{t^{2k}\over(k!)^2}I_{d-1}(W)R^{2(d-1)}+O(R^{2d-3}).$$
\end{corollary}
Especially for the practically relevant case $d=3$, we have for $k=1$ and $k=2$
\begin{eqnarray}
\nonumber \Var(\Sigma_{V_{1};t}^{W_R}) &=& {\pi^2\over 16}t^2I_2(W)R^4+O(R^3),\\
\nonumber \Var(\Sigma_{V_{0};t}^{W_R}) &=& {\pi^2\over 64}t^4I_2(W)R^4+O(R^3).
\end{eqnarray}
In general, $I_{d-1}(W)$ or equivalently $E_2(W)$ is rather difficult to evaluate explicitly. However, for the unit ball $W=B^d$ in ${\Bbb R}^d$ we have by applying \cite{SW}, Theorem 8.6.6 (with a corrected constant) $$I_{d-1}(B^d)={d2^{d-2}}{\kappa_d\kappa_{2d-2}\over\kappa_{d-1}}.$$ For example for $d=3$ this gives us the value $I_2(B^3)=4\pi^2$. For the interpretation of computer simulations it is of particular interest to evaluate the chord power integral $I_{2}(C_a^3)$ for a $3$-dimensional cube $C_a^3$ with side length $a>0$. It can be shown that the numerical value of $I_{2}(C_a^3)$ is given by $5.6337\cdot a^4$.

\section{Second Order Structure of Face Measures}\label{sec2ndMEAS}

 In this section we focus our attention on the spatial pair-correlation structure for the processes of maximal
 polytopes of arbitrary dimensionalities induced by $Y(t),$ arising for $j=0,\ldots,d-1$ as random $j$-volume
 measures concentrated on the union of all $j$-faces of $(d-1)$-dimensional maximal polytopes of $Y(t)$. 
 The nature of the so defined face measures is somewhat different than that of cumulative intrinsic volume
 processes considered above, in particular the face measures keep track not just of the cumulative numeric
 characteristics of STIT tessellations but also of their spatial profile, moreover it should be emphasized
 that in general the total mass of a face measure may be quite unrelated to the cumulative intrinsic volume
 of the corresponding order, as for example $\Sigma_{V_0}(Y(t,W))$, the number of $(d-1)$-dimensional maximal polytopes of
 $Y(t)$ in $W\subset{\Bbb R}^d$, is not deterministically related to the number of vertices of $Y(t,W)$
 as soon as $d > 2$ etc. However, we decided to consider the face measures in this paper as they 
 are of interest in their own right and supplement our results from the other sections, showing 
 the power and versatility of the methods developed in this paper. \\
 For $j=0,1,\ldots,d-1$ we consider the (random) $j$-th order face measure
 $$ {\cal V}_{j;t} := \sum_{e \in \MaxFaces_j(Y(t))} \Vol_j(\cdot \cap e)
   = \sum_{f \in \MaxFacets_{d-1}(Y(t))} v^f_j $$
 where $\Vol_j(\cdot \cap e)$ is the $e$-truncated $j$-volume measure 
 $ [\Vol_j(\cdot \cap e)](A) = \Vol_j(A \cap e),\; A \subseteq {\Bbb R}^d, $
 whereas
 $$ v^f_j := \sum_{e \in \Faces_j(f)} \Vol_j(\cdot \cap e),$$ where by $\Faces_j(f)$ we mean the collection of all $j$-faces of the $(d-1)$-dimensional polytope $f$. We also
 abuse the notation by putting ${\cal V}_{d;t} := \Vol_d(\cdot).$ We shall be interested
 in the covariance measures $\Cov({\cal V}_{i;t},{\cal V}_{j;t}),\; i,j=0,\ldots,d-1$,
 on $({\Bbb R}^d)^2$ given by
 \begin{equation}\label{CovMeasDef}
 \nonumber \langle g \otimes h, \Cov({\cal V}_{i;t},{\cal V}_{j;t}) \rangle = \Cov( \langle g, {\cal V}_{i;t} \rangle,
  \langle h, {\cal V}_{j;t} \rangle)
 \end{equation}
 for all bounded measurable $g,h : {\Bbb R}^d \to {\Bbb R}$ with bounded support, where 
 the standard duality notation $\langle \cdot, \cdot \rangle$ is used for integration 
 $\langle \phi , \mu \rangle = \int \phi d\mu.$ Denote by ${\Bbb Q}^{Y(t)}_{d-1}$ 
 the law of the typical $(d-1)$-dimensional maximal polytope, of the
 STIT tessellation $Y(t)$  and let $\lambda_{Y(t)}^{(d-1)}$ be the corresponding
 facet density. The following Proposition is crucial for this section:
 \begin{proposition}\label{CovRecursionLemma}
  For $i,j \in \{ 0, \ldots, d-1 \}$ we have
  $$ \Cov({\cal V}_{i;t},{\cal V}_{j;t}) = \lambda^{(d-1)}_{Y(t)} \int_{{\Bbb R}^d}  \int
  v^{(f+x)}_i \otimes v^{(f+x)}_j d{\Bbb Q}^{Y(t)}_{d-1}(df)
      dx  + $$
 \begin{equation}\label{CovMeasRecursion}
     \int_0^t [\gamma_{i+1} \Cov({\cal V}_{i+1;s},{\cal V}_{j;s}) + 
                   \gamma_{j+1} \Cov({\cal V}_{i;s},{\cal V}_{j+1;s})] ds
  \end{equation}
   with $\gamma_{i+1}$ and $\gamma_{j+1}$ given as in (\ref{GAMMAJ}). 
  \end{proposition}

 \paragraph{Proof:}
 For bounded measurable and boundedly supported $g$ and $h,$ and for
 $i,j \in \{0,\ldots,d-1\}$  consider the facet functionals 
 $$ J^g_i(f) = \langle g, v^f_i \rangle = 
      \sum_{e \in \Faces_i(f)} \langle g, \Vol_i(\cdot \cap e) \rangle
     = \sum_{e \in \Faces_i(f)} \int_e g(x) \Vol_i(dx) $$ and $J^h_j(\cdot)$
 defined analogously. Then, choosing some compact convex $W$ containing the
 supports of $g$ and $h$ in its interior we see that, recalling (\ref{SIGMADEF}) and (\ref{ADEF}), 
 $$ A_{J^g_i}(Y) = \int_{[W]} \sum_{f \in \Cells(Y \cap H)} J^g_i(f) \Lambda(dH) = 
     \sum_{e \in \MaxFaces_{i+1}(Y)} \left\langle g,  \int_{[W]} \Vol_i(\cdot \cap (e \cap H)) \Lambda(dH) 
     \right\rangle  $$
 \begin{equation}\label{AforMeas}
  = \gamma_{i+1} \sum_{e \in \MaxFaces_{i+1}(Y)} \langle g, \Vol_{i+1}(\cdot \cap e) \rangle 
  = \gamma_{i+1} \Sigma_{J^g_{i+1}}(Y),
 \end{equation}
 where we have used the Crofton formula (\ref{EQCROFTON}). Applying now (\ref{MART2}) with $\phi = J^g_i$ and $\psi = J^h_j,$ we get
 upon taking expectations
 $$
  \Cov(\Sigma_{J^g_i}(Y(t,W)),\Sigma_{J^h_j}(Y(t,W))) =
  \int_0^t {\Bbb E} A_{J^g_i J^h_j}(Y(s,W)) ds + $$ $$
  \int_0^t [\Cov(A_{J^g_i}(Y(s,W)),\Sigma_{J^h_j}(Y(s,W))) + \Cov(\Sigma_{J^g_i}(Y(s,W)),A_{J^h_j}(Y(s,W)))] ds
 $$
 and thus, in view of (\ref{AforMeas}),
 $$ \Cov(\Sigma_{J^g_i}(Y(t,W)),\Sigma_{J^h_j}(Y(t,W))) =
     \int_0^t {\Bbb E} A_{J^g_i J^h_j}(Y(s,W)) ds + $$ $$
     \int_0^t [\gamma_{i+1} \Cov(\Sigma_{J^g_{i+1}}(Y(s,W)),\Sigma_{J^h_j}(Y(s,W))) 
               + \gamma_{j+1} \Cov(\Sigma_{J^g_i}(Y(s,W)),\Sigma_{J^h_{j+1}}(Y(s,W)))] ds. $$
 Using that $\Sigma_{J^g_i}(Y(t,W)) = \langle g, {\cal V}_{i;t} \rangle$ and similar relationships,
 we end up with
$$
  \langle g \otimes h, \Cov({\cal V}_{i;t},{\cal V}_{j;t}) \rangle 
  = \int_0^t {\Bbb E} A_{J^g_i J^h_j}(Y(s,W)) ds + $$
  \begin{equation}\label{CovAux1}
     \left\langle g \otimes h, \int_0^t [\gamma_{i+1} \Cov({\cal V}_{i+1;s},{\cal V}_{j;s}) + \gamma_{j+1}
     \Cov({\cal V}_{i;s},{\cal V}_{j+1;s})]  ds \right\rangle.
 \end{equation}
 Thus, to establish (\ref{CovMeasRecursion}) it is enough to show that
 \begin{equation}\label{CovAux2}
  \int_0^t {\Bbb E} A_{J^g_i J^h_j}(Y(s,W)) ds = \left\langle g \otimes h, 
  \lambda^{(d-1)}_{Y(t)} \int_{{\Bbb R}^d} \int
  v^{(f+x)}_i \otimes v^{(f+x)}_j {\Bbb Q}^{Y(t)}_{d-1}(df) dx  \right\rangle.
\end{equation} 
 To prove (\ref{CovAux2}) use (\ref{STITSECTION}) and  write, recalling (\ref{ADEF}), 
 \begin{eqnarray}
 \nonumber {\Bbb E}A_{J^g_i J^h_j}(Y(s,W)) &=& \int_{[W]} {\Bbb E} \sum_{f \in \Cells(Y(s,W) \cap H)} J^g_i J^h_j(f)\Lambda(dH)\\
 &=& \int_{[W]} \int J^g_i J^h_j(f) {\Bbb M}^{Y(s,W) \cap H}(df) \Lambda(dH),\label{EQXXX}
 \end{eqnarray}
 where 
 ${\Bbb M}^{Y(s,W) \cap H}$ is the cell intensity measure for the sectional STIT tessellation $Y(s,W) \cap H$, compare with (\ref{CellIntensityMeasure}). Using Proposition \ref{PropXX} (a) we are led to $$\int_{[W]} \int J^g_i J^h_j(f) {\Bbb M}^{Y(s,W) \cap H}(df) \Lambda(dH)=\int_{[W]} \int J^g_i J^h_j(f) {\Bbb M}^{\PHT(s,W) \cap H}(df) \Lambda(dH).$$ Now, using Slivnyak's Theorem (cf. \cite[Thm. 3.2.5]{SW}) we obtain for any bounded measurable function $\phi$ on the space of $(d-1)$-dimensional polytopes
\begin{eqnarray}
\nonumber \int_{[W]}\int\phi(f){\Bbb M}^{\PHT(1,W)\cap H}(df)\Lambda(dH) &=& \int\int_{[c]}\phi(c\cap H)\Lambda(dH){\Bbb M}^{\PHT(1,W)}(dc)\\
\nonumber &=& \int\phi(f){\Bbb F}_{d-1}^{\PHT(1,W)},
\end{eqnarray}
where ${\Bbb M}^{\PHT(1,W)\cap H}$ is the cell intensity measure of $\PHT(1,W)\cap H$, whereas ${\Bbb M}^{\PHT(1,W)}$ is the cell intensity measure and ${\Bbb F}^{\PHT(1,W)}_{d-1}$ is the $(d-1)$-face intensity measure of a Poisson hyperplane tessellation $\PHT(1,W)$ within $W$ having intensity measure $\Lambda$. Thus, replacing $\Lambda$ by $s\Lambda$, $0<s\leq t$, we arrive at $$\int_{[W]}\int\phi(f){\Bbb M}^{\PHT(s,W)\cap H}(df)\Lambda(dH)={1\over s}\int\phi(f){\Bbb F}_{d-1}^{\PHT(s,W)}(df),$$ whence continuing (\ref{EQXXX}) we find
\begin{equation}
 {\Bbb E}A_{J^g_i J^h_j}(Y(s,W)) = \frac{1}{s} \int J^g_i J^h_j(f) {\Bbb F}^{\PHT(s,W)}_{d-1}(df). \label{AAux1}
\end{equation}
 Hence, by applying Proposition \ref{PropXX} (b) we get from (\ref{AAux1}),
 \begin{equation}\label{AAux2}
  \int_0^t {\Bbb E} A_{J^g_i J^h_j}(Y(s,W)) ds = \int J^g_i J^h_j(f) {\Bbb F}^{Y(t,W)}_{d-1}(df), 
 \end{equation}
 where, recall, ${\Bbb F}^{Y(t,W)}_{d-1}$ is the $(d-1)$-dimensional maximal polytope intensity measure for $Y(t,W)$, see (\ref{FaceIntensityMeasure}). In view of (\ref{AAux2}) 
 the required relation (\ref{CovAux2}) follows now directly
 by the definition of typical $(d-1)$-dimensional maximal polytope \cite[Chap. 4.1, Chap. 10]{SW} and the definition of $J^g_i,$ 
 upon taking into account that the supports of both $g$ and $h$ 
 are contained in the interior of $W.$ Putting (\ref{CovAux1}) and (\ref{CovAux2})
 together yields
 (\ref{CovMeasRecursion}) and thus completes the proof. \hfill $\Box$ \\ \\
 To proceed note that ${\cal V}_{d;s}$ is a constant measure and therefore any covariances 
 involving it vanish. Thus, using (\ref{CovMeasRecursion}) in Proposition \ref{CovRecursionLemma}
 and arguing as in the derivation of Theorem \ref{thmconvexact} from the crucial relation
 (\ref{PrincipalCov}) we obtain
 \begin{theorem}\label{thmCovMeas}
  For $k,l \in \{0,\ldots,d-1\}$ we have
  $$ \Cov({\cal V}_{d-1-k;t},{\cal V}_{d-1-l;t}) = 
       \sum_{m=0}^k \sum_{n=0}^l {{k+l-m-n} \choose {k-m}}
       \left(\prod_{i=m+1}^k \gamma_{d-i} \right) \left( \prod_{j=n+1}^l \gamma_{d-j} \right) $$
  \begin{equation}\label{CovMeasFullExpr}
    \times\int_{{\Bbb R}^d} \int [v^{(f+x)}_{d-1-m} \otimes v^{(f+x)}_{d-1-n}] 
       {\cal I}^{k+l-m-n}\left(\lambda^{(d-1)}_{Y(\cdot)} {\Bbb Q}^{Y(\cdot)}_{d-1};t \right)(df) dx
  \end{equation}
   with the usual convention that ${0\choose 0} = 1$ and $\prod_{k+1}^k \ldots = \prod_{l+1}^l \ldots \equiv 1.$
  \end{theorem} 
  It is interesting to note that for the particular case $k=l=0$ we get
 \begin{corollary}\label{CovCor1}
  We have
  $$ \Cov({\cal V}_{d-1;t}) = \lambda^{(d-1)}_{Y(t)} 
       \int_{{\Bbb R}^d} \int [v_{d-1}^{(f+x)}\otimes v_{d-1}^{(f+x)}]{\Bbb Q}^{Y(t)}_{d-1}(df) dx. $$
 \end{corollary}
  This means that for $k=l=0$ the covariance measure of the surface area process ${\cal V}_{d-1;t}$ 
  coincides with that of the surface area process induced by the homogeneous and isotropic Boolean
  model with grain distribution ${\Bbb Q}^{Y(t)}_{d-1}$ and with grain density $\lambda^{(d-1)}_{Y(t)},$
  a result first established by Weiss, Ohser and Nagel \cite{NOW} in the special planar case by completely different methods. Recall further from Theorem 3 in \cite{ST} that $${\Bbb Q}_{d-1}^{Y(t)}=\int_0^t{ds^{d-1}\over t^d}{\Bbb Q}_{d-1}^{\PHT(s)}ds,$$ where ${\Bbb Q}_{d-1}^{\PHT(s)}$ stands for the distribution of the typical facet of the Poisson hyperplane tessellation $\PHT(s)$ with surface intensity $s$ (this is to say, the hyperplane process has intensity measure $s\Lambda$). But this random polytope has the same distribution ${\Bbb Q}_{\TypicalCell}^{\PHT(\gamma_{d-1}s,{\Bbb R}^{d-1})}$ as the typical cell of a Poisson hyperplane tessellation in ${\Bbb R}^{d-1}$ with surface intensity $\gamma_{d-1}s$ with $\gamma_{d-1}$ given by (\ref{GAMMAJ}). Thus, $$\Cov({\cal V}_{d-1;t})=\lambda_{Y(t)}^{(d-1)}\int_0^t{ds^{d-1}\over t^d}\int_{{\Bbb R}^d}\int[v_{d-1}^{c+x}\otimes v_{d-1}^{c+x}]{\Bbb Q}_{\TypicalCell}^{\PHT(\gamma_{d-1}s,{\Bbb R}^{d-1})}(dc)dxds$$ and the covariance measure is reduced to known quantities.\\ For $k+l > 0$ the
  situation becomes more complicated and mixtures of typical cell distributions corresponding to different
  time moments arise in the right hand side of (\ref{CovMeasFullExpr}). As in Section \ref{sec2ndORD} above,
  also here explicit calculations are precluded
  for $k+l > 0$ due to the lack of known formulae for mixed moments of general order intrinsic volumes
  of Poisson cells.

\section{Central Limit Theory}\label{secCLT}

 In this section we present a central limit theory for the suitably rescaled intrinsic volume processes
 $\Sigma_{V_i;t}.$ In strong contrast to the two-dimensional case considered in \cite{ST2}, where
 a classical Gaussian limit behaviour is observed, for the case $d \geq 3$ in focus of this paper the
 situation is very different and non-Gaussian limits arise. In this context, to proceed with a full discussion
 below, we recapitulate first some facts already known from \cite{ST} in a way specialized for our
 present purposes.  Define the rescaled intrinsic volume processes $({\cal S}^{R,W}_{V_i;t})_{t\in[0,1]}$ for $i=0,\ldots,d-1$ by
 \begin{equation}\label{RIV}
 \nonumber {\cal S}^{R,W}_{V_i;t} = R^{-(d-1)} \bar\Sigma_{V_i;t+\log R/R}^{W_R} = 
  R^{-(d-1)} \bar \Sigma_{V_i}(Y(t+\log R/R,W_R)),\; t \in [0,1].
 \end{equation}
 Note that the shift by $\log R/R$ in time argument here is of technical importance as placing the time origin
 just after the very initial {\it big bang} phase of the MNW-construction, where the dominating fluctuations
 of $\bar\Sigma_{V_{d-1};(\cdot)}$ arise, but where nothing of asymptotic significance happens for
 $\bar\Sigma_{V_{i};(\cdot)},\; i < d-1.$ Since the big bang phase evolution for $\bar\Sigma_{V_{d-1};(\cdot)}$ 
 has been considered in full detail in \cite[Sec. 5]{ST}, in this paper we only focus on the later phase 
 $[\log R/R,1].$ Putting together the present Section with \cite[Sec. 5]{ST} reveals remarkable
 richness of the complete asymptotic picture for STIT tessellations in large windows.\\
 We begin with the observation that follows by the theory developed in Subsection 5.3 of \cite{ST}:
 \begin{proposition}\label{VD1}
  The process  $({\cal S}^{R,W}_{V_{d-1};t})_{t \in [0,1]}$ converges in law, as $R \to \infty,$ in the
  space ${\cal D}[0,1]$ of right-continuous functions with left-hand limits (c\`adl\`ag) on $[0,1]$ 
  endowed with the usual Skorokhod $J_1$-topology \cite[Chap. VI.1]{JS} to the constant process
  $t \mapsto \xi$ where $\xi := \Xi(W)$ is a certain non-Gaussian square-integrable random variable
  with variance $\frac{d-1}{2} E_2(W) = \frac{1}{d-2} I_{d-1}(W).$
 \end{proposition} 
 \paragraph{Proof:}
  The relation (85) in \cite{ST} implies that ${\cal S}^{R,W}_{V_{d-1};1-\log R/R}$ converges in
  law to $\xi$ as $R \to \infty.$ On the other hand, recall that ${\cal S}^{R,W}_{V_{d-1};t}$ is
  a martingale in view of (\ref{MartExpl1ord}). Moreover, by (\ref{TOTALVAR}) and the definition of ${\cal S}^{R,W}_{V_{d-1};t}$ combined with the scaling property of STIT tessellations we have
  \begin{eqnarray}
  \nonumber \Var({\cal S}^{R,W}_{V_{d-1};t}) &=& \Var (R^{-(d-1)}\bar\Sigma_{V_{d-1}}(Y(t+\log R/R,W_R)))\\
  \nonumber &=& \Var(\bar \Sigma_{V_{d-1}}(Y(R(t+\log R/R),W)))\\
  \nonumber    &=& \frac{d-1}{2} \int_W \int_W \frac{1-e^{-\frac{2\kappa_{d-1}}{d\kappa_d} [Rt + \log R] ||x-y||}}{||x-y||^2}
       dx dy
  \end{eqnarray}
  and hence $\Var({\cal S}_{V_{d-1};1}^{R,W}) - \Var({\cal S}_{V_{d-1};0}^{R,W})$ tends to $0$ as
  $R \to \infty.$ Consequently, the asymptotic constancy of the limit process follows now by 
  Doob's $L^2$-maximal inequality \cite[Thm 3.8(iv)]{KS} which completes our argument. \hfill $\Box$\\ \\

\begin{remark}\label{NonGaussExpl}
The convergence in law and non-Gaussianity results referred to above come from the paper \cite{ST}. In order to keep the present work formally self-contained we give here a very brief sketch of the main arguments.\\ The convergence is readily guaranteed by the 
 martingale convergence theory as ibidem. The crucial point is the non-Gaussianity. Here,
 the idea, as described in detail in Section 5.3 of \cite{ST}, relies on showing that 
 $\xi$ has its tails much heavier than normal random variables. We proceed, roughly speaking, 
 by constructing a suitable class of {\it initial frame tessellations} of the body $W$ 
 with the properties that
 \begin{itemize}
  \item the number of frame hyperplanes hitting $W$ equals $N\in{\Bbb N}$,
  \item the frame hyperplanes intersect within $W$ only very seldom, which implies that $W$
          gets subdivided into $\Theta(N)$ cells (recall the Landau notation).
 \end{itemize}
 This can be achieved by choosing $d$ principal directions and keeping all hyperplanes in the collection
 under construction approximately aligned to these directions. Importantly, 
 it turns out that this way we can get
 the lower bound $\exp(-O(N \log N))$ for the probability that the real initial frame arising in
 the very initial phase for $Y(\cdot,W)$ -- usually called {\it big bang}, see below -- upon suitable
 spatial re-scaling does fall into the desired class. Now, denoting 
 the last event by ${\cal E}_N$, we see that on ${\cal E}_N$ the total deviation 
 $\xi \approx {\cal S}^{R,W}_{V_{d-1};t}$ gets decomposed into
 \begin{itemize}
  \item a sum of independent contributions $\xi_i$ coming from all respective cells $W_i$
          of the frame tessellation,
  \item a correction term of order $\Theta(N)$ due to the centred joint contribution of all hyperplanes
          of the frame tessellation (observe that each such hyperplane necessarily contributes
          $\Theta(1)$, whereas ${\Bbb E}\xi = \Theta(1)$).
 \end{itemize}
 Finally we can use the knowledge of $\Var(\xi)$ and $\Var(\xi_i)$ and properly exploit the
 independence of the random variables $\xi_i$, to conclude that on ${\cal E}_N$ we can find a further sub-event
 ${\cal E}'_N \subseteq {\cal E}_N$ with ${\Bbb P}({\cal E}'_N) = \Theta({\Bbb P}({\cal E}_N))
 = \exp(-O(N \log N))$ such that $\{ \xi \geq N \} \supseteq {\cal E}'_N$ and hence
 ${\Bbb P}(\xi > N) = \exp(-O(N \log N))$, which could not hold if $\xi$ were Gaussian.
 We refer the reader to \cite{ST} for further details.
\end{remark}
  Having characterized the asymptotic behaviour of $({\cal S}^{R,W}_{V_{d-1};t})_{t \in [0,1]}$ we
  are now prepared to describe the full joint asymptotics of all intrinsic volume processes, which is the
  main result of this section.
  \begin{theorem}\label{IVLimTheo}
   The vector $({\cal S}_{V_{d-1};t}^{R,W},{\cal S}_{V_{d-2};t}^{R,W},\ldots,
   {\cal S}_{V_0;t}^{R,W})_{t \in [0,1]}$ converges in law, as $R \to \infty,$ in the
   space ${\cal D}([0,1];{\Bbb R}^d)$ of ${\Bbb R}^d$-valued c\`adl\`ag functions on
   $[0,1]$ endowed with the usual Skorokhod $J_1$-topology \cite[Chap.VI.1]{JS}, to
   the stochastic process 
   \begin{equation}\label{LimProc}
    t \mapsto \left(\xi,\gamma_{d-1} t \xi, \frac{\gamma_{d-1} \gamma_{d-2} t^2}{2!} \xi, 
   \frac{\gamma_{d-1} \gamma_{d-2} \gamma_{d-3} t^3}{3!} \xi, \ldots, 
   \frac{\gamma_{d-1} \ldots \gamma_1 t^{d-1}}{(d-1)!} \xi\right),\;\; t \in [0,1]. 
  \end{equation}
 \end{theorem}
  Before proceeding with the proof we discuss certain striking features of the phenomenon 
  described in Theorem \ref{IVLimTheo}. Namely, although all intrinsic volume processes
  ${\cal S}^{R,W}_{V_i;t}$ exhibit fluctuations of order $R^{d-1},$ the nature of these
  fluctuations differs very much between $i=d-1$ and $i < d-1.$ 
  \begin{itemize}
   \item As shown in Subsection 5.3 in \cite{ST} and in Proposition \ref{VD1} above, the
     leading-order deviations of $\bar \Sigma_{V_{d-1};t}^{R,W}$ in large $R$ asymptotics
     arise very early in the course of the MNW-construction, in its initial stages usually referred
     to as the {\it big bang phase}. Here, this is the time period $[0,\log R/R].$ During the later
     stages of the construction,  i.e. say the time interval $(\log R/R,1],$ the variance increase is of
     lower order and any newly arising fluctuations are negligible compared to those originating from
     the {\it big bang}. In the asymptotic picture this means that the initial fluctuation remains {\it frozen}
     throughout the rest of the dynamics, whence the constant limit for
     $({\cal S}^{R,W}_{V_{d-1};t})_{t \in [0,1]}$  as $R \to \infty.$  
  \item In contrast, the leading-order deviations of $\bar \Sigma_{V_i;t}^{R,W}$ for $i < d-1$ arise and cumulate
           constantly in $t$ with deterministic polynomial rates depending on $i,$ times the initial
           {\it big bang} fluctuation of the process $\bar \Sigma_{V_{d-1};(\cdot)}^{R,W}$ which,
           in this sense, {\it stores the entire randomness} of the intrinsic volume vector. The mechanism
           determining the dependence of fluctuations of intrinsic volume processes of orders $i < d-1$ 
           given those for $d-1$ and the resulting form of the limit process  (\ref{LimProc}) will be discussed
           in the sequel. 
 \end{itemize}
             
 \paragraph{Proof of Theorem \ref{IVLimTheo}:}
  The crucial step of the proof relies on considering for each $j=0,\ldots,d-1$ the auxiliary process 
 \begin{equation}\label{SigmaAux}
 \hat\Sigma_{V_j;t} := \bar \Sigma_{V_j;t} - \int_0^t \bar A_{V_j;s} ds,
 \end{equation}
 which is a centred $\Im_t$-martingale by (\ref{MART1}) and which is the same as
 \begin{equation}\label{SigmaAux2}
 \hat\Sigma_{V_j;t} = \bar \Sigma_{V_j;t} - \gamma_{j+1} \int_0^t \bar \Sigma_{V_{j+1};s} ds 
 \end{equation}
 by (\ref{BarAj}). The idea is to show that for $j < d-1$ the processes
 $\hat\Sigma_{V_j;t}$ are of negligible order in large $R$ asymptotics. Indeed, upon squaring
 and taking expectations we get
 \begin{equation}\label{VarSigmaAux}
  {\Bbb E}(\hat\Sigma_{V_j;t})^2 = \Var(\Sigma_{V_j;t})  - 2 \gamma_{j+1} \int_0^t 
  {\Bbb E} \bar \Sigma_{V_j;t} \bar \Sigma_{V_{j+1};s} ds + \gamma_{j+1}^2 \int_0^t \int_0^t 
  {\Bbb E} \bar \Sigma_{V_{j+1};s} \bar  \Sigma_{V_{j+1};u} du ds.
 \end{equation}
 Using that
 $$ {\Bbb E}(\bar\Sigma_{V_{j};t}|\Im_s) = \bar\Sigma_{V_j;s} + \gamma_{j+1} \int_s^t 
     {\Bbb E}(\bar\Sigma_{V_{j+1};u}|\Im_s) du, $$
 as follows by the martingale property of $\hat\Sigma_{V_j;t},$ and inserting this to (\ref{VarSigmaAux})
 we are led to
 $$ {\Bbb E}(\hat\Sigma_{V_j;t})^2 = \Var(\Sigma_{V_j;t}) - 2 \gamma_{j+1} \int_0^t
  {\Bbb E} \bar\Sigma_{V_j;s} \bar\Sigma_{V_{j+1};s} ds - $$
 $$ 2 \gamma_{j+1}^2 \int_0^t  \int_s^t {\Bbb E} ({\Bbb E}(\bar\Sigma_{V_{j+1};u}|\Im_s))  
        \bar\Sigma_{V_{j+1};s} du ds +  
       \gamma_{j+1}^2 \int_0^t \int_0^t  {\Bbb E} \bar \Sigma_{V_{j+1};s} \bar  \Sigma_{V_{j+1};u} du ds
 $$
 with the last two terms cancelling. Thus,
 $$ {\Bbb E} (\hat\Sigma_{V_j;t})^2 = \Var(\Sigma_{V_j;t}) - 2 \gamma_{j+1} \int_0^t \Cov(\Sigma_{V_j;s},
     \Sigma_{V_{j+1};s}) ds, $$
 whence, by (\ref{PrincipalCov}), 
 \begin{equation}\label{VarSigmaAux2}
  {\Bbb E} (\hat\Sigma_{V_j;t})^2 = \int_0^t {\Bbb E} A_{V_j^2;s} ds.
 \end{equation}
 In view of Proposition \ref{PropAux1} for $j < d-1$ we have
 $$ \int_0^t {\Bbb E} A^{W_R}_{V_j^2;s} ds = o(R^{2(d-1)}),\; t \in [0,1] $$
 and thus the relation (\ref{VarSigmaAux2}) implies
 $$ \Var(\hat\Sigma_{V_j;t}^{W_R}) = o(R^{2(d-1)}),\; t \in [0,1] $$
 for $j < d-1.$ Hence finally
 \begin{equation}\label{FinalNegl}
   {\Bbb E} \sup_{t \in [0,1]}  (\hat\Sigma_{V_j;t}^{W_R})^2 = o(R^{2(d-1)})
 \end{equation}
 for $j < d-1$ by Doob's $L^2$-maximal inequality \cite[Thm 3.8(iv)]{KS} applied
 to the martingale $\hat\Sigma^{W_R}_{V_j;t}$ defined as in (\ref{SigmaAux}, \ref{SigmaAux2})
 with $W$ replaced by $W_R$ there according to our usual notational convention.\\
 With (\ref{FinalNegl}) it is now easy to complete the proof. Indeed, since the normalization in
 the definition of intrinsic volume processes ${\cal S}_{V_j;t}^{R,W}$ involves a prefactor $R^{-(d-1)},$ 
 the relation (\ref{FinalNegl}) allows us to recursively substitute
 $$ \gamma_{j+1} \int_0^t {\cal S}_{V_{j+1};s}^{R,W} ds $$ for
 $ {\cal S}_{V_j;t}^{R,W} $, as soon as $j < d-1$, without affecting the large $R \to \infty$
 asymptotics in law (note that the technically motivated shift by $\log R/R \to 0$ in time argument of
 the rescaled process ${\cal S}_{V_j;t}^{R,W}$ is asymptotically negligible as inducing only a
 negligible $L^2$-difference precisely calculated in Theorem \ref{thmconvexact}).
  The application of such recursive substitutions combined with Proposition \ref{VD1}
 completes the proof of Theorem \ref{IVLimTheo}.  \hfill $\Box$\\ \\
 It is crucial to emphasize at this point that, in intuitive terms, under the normalization of
 ${\cal S}_{V_j;t}^{R,W},$ $j=0,\ldots,d-1,$  the mechanism governing the rise of fluctuations
 of intrinsic volume processes of orders $j<d-1$ given those for $V_{d-1}$ reduces effectively
 to the simple approximation 
 $$ {\cal S}_{V_j;t}^{R,W} \approx \gamma_{j+1} \int_0^t {\cal S}_{V_{j+1};s}^{R,W} ds $$
 and its recursive application. Of course this {\it simple approximation} follows itself by rather
 non-trivial arguments above. Clearly, this description only characterizes the leading order
 fluctuations as considered in Theorem \ref{IVLimTheo} although our tools should allow
 a more delicate characterization of non-leading lower order fluctuations
 as well, see e.g. Subsection 5.1 in \cite{ST} for the particular case of $V_{d-1}.$ 
 
\subsection*{Acknowledgement} 

We are grateful to Claudia Redenbach and Joachim Ohser who provided the two pictures of STIT tessellations. We would also like to thank the anonymous referees for their comments and suggestions as well as the Editor-in-Chief for his help.\\ The first author was supported by the Polish Minister of Science and Higher Education grant N N201 385234 (2008-2010). The second author was supported by the Swiss National Science Foundation, Grant PP002-114715/1.


\begin{thebibliography}{30}\small
\bibitem{Bertoin}
J. Bertoin, Random Fragmentation and Coagulation Processes, Cambridge University Press (2006).

\bibitem{Cowan}
R. Cowan, New classes of random tessellations arising from iterative division of cells, Adv. Appl. Probab. \textbf{42} (2010) 26--47.

\bibitem{JS}
J. Jacod, A.N. Shiryaev, Limit theorems for Stochastic Processes, second ed., Springer (2003). 

\bibitem{Kallenberg}
O. Kallenberg, Foundations of Modern Probability, second ed., Springer (2003).

\bibitem{KS}
I. Karatzas, S.E. Shreve, Brownian Motion and Stochastic Calculus, second ed., Springer (1998).

\bibitem{MNW}
J. Mecke, W. Nagel, W. Weiss, A global construction of homogeneous random planar tessellations that are stable under iteration, Stochastics \textbf{80} (2008) 51--67.

\bibitem{Miles}
R.E. Miles, Poisson flats in Euclidean spaces. Part II: Homogeneous Poisson flats and the complementary theorem, Adv. Appl. Probab. \textbf{3} (1971) 1--43.

\bibitem{NW05}
W. Nagel, V. Weiss, Crack STIT tessellations: characterization of stationary random tessellations stable with respect to iteration, Adv. Appl. Probab. \textbf{37} (2005) 859--883.

\bibitem{NW06}
W. Nagel, V. Weiss, STIT tessellations in the plane, Rendiconti del circulo matematico di Palermo, Serie II, Suppl. 77 (2006) 441--458.

\bibitem{NMOW}
W. Nagel, J. Mecke, J. Ohser, V. Weiss, A tessellation model for crack pattern on surfaces, Image. Anal. Stereol. \textbf{27} (2008) 73--78.

\bibitem{SW}
R. Schneider, W. Weil, Stochastic and Integral Geometry, Springer (2008).

\bibitem{ST}
T. Schreiber, C. Th\"ale, Typical geometry, second-order properties and central limit theory for iteration stable tessellations, arXiv: 1001.0990 [math.PR] (2010).

\bibitem{ST2}
T. Schreiber, C. Th\"ale, Second-order properties and central limit theory for the vertex process of iteration infinitely divisible and iteration stable random tessellations in the plane, Adv. Appl. Probab. \textbf{42}, 913--935 (2010).

\bibitem{NOW}
V. Weiss, J. Ohser, W. Nagel, Second moment measure and K-function for planar STIT tessellations, Image Anal. Stereol. \textbf{29} (2010) 121--131.

\end{thebibliography}
\end{document}